\documentclass[a4paper,12pt]{article}

\usepackage[latin1]{inputenc}
\usepackage{amsmath}
\usepackage{amsfonts}
\usepackage{amssymb}
\usepackage{amsthm}
\usepackage{fontenc}
\usepackage[pdftex]{graphicx}
\usepackage[left=2.5cm,right=2.5cm,top=2cm,bottom=2cm]{geometry}
\usepackage{cite}
\usepackage{makeidx}
\usepackage{graphicx}
\usepackage{tikz}
\usepackage{pgf}

\newcommand{\diam}{\mathrm{diam}}
\newcommand{\dist}{\mathrm{dist}}

\newtheorem{definition}{Definition}
\newtheorem{theorem}{Theorem}
\newtheorem{lemma}{Lemma}

\newtheorem{further}{Further results}

\begin{document}

\title{Marcinkiewicz Exponent and Boundary Value Problems in Fractal Domains of $\mathbb{R}^{n+1}$}

\author{Carlos Daniel Tamayo Castro}
\date { \small{Instituto de Matem\'aticas. Universidad Nacional Aut\'onoma de M\'exico, Ciudad M\'exico, M\'exico.\\
	cdtamayoc@comunidad.unam.mx}
}

\maketitle

\begin{abstract} 
This paper aims to study the jump problem for monogenic functions in fractal hypersurfaces of Euclidean spaces. The notion of the Marcinkiewicz exponent has been taken into consideration. A new solvability condition is obtained, basing the work on specific properties of the Teodorescu transform in Clifford analysis. It is shown that this condition improves those involving the Minkowski dimension.
\vspace{0.3cm}

\small{
 }  
\end{abstract}

\section{Introduction}\label{sec1}

The Riemann boundary value problem in Complex Analysis is widely used in many branches of Mathematics and Physics. The classical references here are \cite{Gajov, Lu, Mu}. In the solution to the Riemann boundary value problem, the Cauchy type integral    
\begin{equation}\label{Cau_Int_Eq_CA}
	(\mathcal{C}_{\gamma}u)(z) = \dfrac{1}{2\pi i}\int\limits_{\gamma}\dfrac{u(\tau)}{\tau - z}d\tau,
\end{equation}
is used as the main mathematical apparatus. It is well known that for every H\"older continuous function $u$ with exponent $\nu > \frac{1}{2}$, $(\mathcal{C}_{\gamma}u)(z)$ has continuous limit values on a rectifiable closed Jordan curve $\gamma$; hence the jump problem is solvable. For a thorough description of old and recent results concerning the geometric conditions on a Jordan curve in the plane that imply the boundedness of the Cauchy type integral boundary behavior, the reader is referred to \cite{ABK2015}.\\ 
However, in the context of non-rectifiable curves, the Cauchy type integral (\ref{Cau_Int_Eq_CA}) has no sense. In contrast, the Riemann boundary value problem is still completely valid. In the early eighties, a complete treatment of this topic was given by B. A. Kats in \cite{BKats83}. There is shown that under the condition $\nu > \frac{\overline{\dim}_{M}(\gamma)}{2}$, the Riemann problem is solvable. Here $\overline{\dim}_{M}(\gamma)$ is the upper Minkowski dimension of the curve $\gamma$, and $\nu$ is the H\"older exponent of a function $u$ associated with the problem. \\
In \cite{1DKats16, 2DKats16}, the Marcinkiewicz exponent is introduced. Using this, there is obtained, once more in the context of Complex Analysis, a solvability condition that improves the conditions mentioned above. \\    
Clifford Analysis provides a natural generalization of Complex Analysis to higher dimensions. One of its main approaches studies functions defined at $\mathbb{R}^{n + 1}$ and valued in the Clifford Algebras $\mathcal{C}\ell(n)$, mainly those that nullify the cliffordian Cauchy-Riemann operator. We follow \cite{TAB21}, calling this approach paravectorial Clifford analysis. On the other hand, we will call vectorial Clifford analysis the approach that studies functions defined at  
$\mathbb{R}^{n}$ with values in $\mathcal{C}\ell(n)$ that nullify the Dirac operator. The reader is referred  to  \cite{BDS, GHS08,MMitrea1994} for a standard account of the theory.\\
Significant obstacles exist when giving a complete treatment to the Riemann boundary value problem for monogenic functions, as seen in \cite[pp 22, 24]{Ryan96}. These  are a direct consequence of the fact that the product of two monogenic functions is not necessarily a monogenic function due to the lack of commutativity in the cliffordian product. This fact explains why an explicit solution to the Riemann boundary value problems has been found only for the jump problem and some slight modifications, where the problematic of the cliffordian multiplication can be essentially avoided, see \cite{AB2001,Brs2001}. In this sense, it is worth pointing out the recent article \cite{TAB21}. There is shown how, in the vectorial approach, the non-commutative product induces substantial differences in the number of solutions. An example of a homogeneous Riemann boundary value problem with constant coefficients is provided with an infinite number of linearly independent solutions.\\
In the context of Clifford Analysis, the research about the solvability of the Riemann problem over fractal domains is focused on the particular case of the jump problem. In this sense, the current results on the plane involving the Minkowski dimension have been generalized to higher dimensions in \cite{APB2007,ABT2007,ABT2006}. Besides, conditions exist that involve the approximate dimension and the $d$-summability; see the article \cite{ABK2013} and its references. Unfortunately, the condition of the $d$-summability does not improve the previous results. However, it has its advantages at the moment of dealing with it. At the same time, the approximate dimension is rather complicated in computations. \\
Nevertheless, to the best of the author's knowledge, there is no research leading to study, in higher dimensions, the relations between the Marcinkiewicz exponent and the jump problem in domains with a fractal boundary.\\       
The main goal of this paper is precisely this, to obtain solvability conditions involving the Marcinkiewicz exponent for the jump problem in Clifford Analysis. Besides, to show that these conditions improve those involving the Minkowski dimension. Also, an example in dimension three is constructed, illustrating that for every value of the Minkowski dimension, there exists a non-numerable class of surfaces where the inequality relating the Marcinkiewicz exponent and the Minkowski dimension is strict. \\
This paper is organized as follows: Section 2 outlines Clifford algebras basic principles and properties, monogenic functions, fractal dimensions, and Teodorescu transform. In Section \ref{Marcinkiewicz}, the Marcinkiewicz exponent in $\mathbb{R}^{n+1}$ is defined. There is proved a lemma that is essential to extend to higher dimensions the inequality involving the Minkowski dimension and the Marcinkiewicz exponent, which is also proved there. Section \ref{Example} is devoted to constructing a class of surfaces. One of the key results of this paper is to show, using these surfaces, that for any possible value of the Minkowski dimension between two and three, there is a non-numerable amount of surfaces where the inequality from Section 3 is exact. The other main achievement is shown in Section \ref{JumpProblem}, where we get conditions for solvability and unicity in a class for the jump problem that improves those conditions involving the Minkowski dimension.      

\section{Preliminaries}\label{Preliminars}
This section briefly presents the results needed to develop this paper.
\subsection{Clifford Algebras and Monogenic Functions}
This subsection has compiled some basic facts concerning Clifford Algebras and Clifford Analysis. For a discussion of this topic, as was mentioned before, the reader is referred to \cite{BDS, GHS08, MMitrea1994}.
\begin{definition}\label{CliffAlg}
	The Clifford algebra associated with $\mathbb{R}^{n}$, endowed with the usual Euclidean metric, is the extension of $\mathbb{R}^{n}$ to a unitary, associative algebra $\mathcal{C}\ell(n)$ over the reals, which is generated, as an algebra, by $\mathbb{R}^{n}$. It is not generated by any proper subspace of $\mathbb{R}^{n}$ and satisfies 
	\begin{eqnarray*}
		x^{2} = -\arrowvert x\arrowvert^{2},   
	\end{eqnarray*}
	for   any  $x \in \mathbb{R}^{n}$.
\end{definition}
Definition \ref{CliffAlg} is equivalent to the next construction. Let $\{e_{j}\}_{j=1}^{n}$ be the standard basis of $\mathbb{R}^{n}$. Multiplication is defined through the basic rule
\begin{equation*}
	e_{i}e_{j} + e_{j}e_{i} = -2\delta_{ij},
\end{equation*}
where $\delta_{ij}$ is the Kronecker delta.

Every $a \in \mathcal{C}\ell(n)$ is of the form $a = \sum\limits_{A \subseteq N}a_{A}e_{A}$, with  $N = \{1,\ldots, n\}$, $a_{A} \in \mathbb{R}$, where $e_{\emptyset} = e_{0} = 1$, $e_{\{j\}} = e_{j}$ and $e_{A} = e_{\beta_{1}}\ldots e_{\beta_{k}}$, for $A = \{\beta_{1}, \ldots \beta_{k}\}$ where $\beta_{k} \in \{1,\ldots, n\}$ and $\beta_{1} < \beta_{2} < \ldots < \beta_{k}$.\\ 
The conjugation is defined by $b \rightarrow \overline{b} := \sum\limits_{A}b_{A}\overline{e_{A}}$  where
\begin{equation*}
	\overline{e_{A}} := (-1)^{k}e_{\beta_{k}}\ldots e_{\beta_{2}}e_{\beta_{1}}.
\end{equation*} 
An algebra norm is defined on $\mathcal{C}\ell(n)$ through $\arrowvert b\arrowvert = \left(\sum_{A}b_{A}^{2}\right)^{\frac{1}{2}}$. With this, $\mathcal{C}\ell(n)$ becomes a Euclidean space. Every point in $\mathbb{R}^{n + 1}$ can be identified with the Clifford number $x = \sum\limits_{i = 0}^{n}x_{i}e_{i}$, named paravector. \\
Functions to be studied will be the ones defined in a domain $G \subset \mathbb{R}^{n + 1}$ valued in $\mathcal{C}\ell(n)$. They have the form
\begin{equation*}
	u(x) = \sum_{A}u_{A}(x)e_{A},
\end{equation*}
where $u_{A}(x)$ are real-valued. From now on, unless the opposite is specified, all functions will be considered Clifford-valued. We say that $u \in C^{k}(G)$ if all the components $u_{A} \in C^{k}(G)$. In general, we say that a Clifford-valued function belongs to a class if all its components belong to that class.\\
The cliffordian Cauchy-Riemann operator in $\mathbb{R}^{n + 1}$ is defined as
\begin{equation*}
	\mathcal{D} := \sum_{i=0}^{n}e_{i}\partial_{i},
\end{equation*}
where $\partial_{i}:= \frac{\partial}{\partial x_{i}}$ is the partial derivative with respect to $x_{i}$. If needed, we shall specify the variable we are applying the operator to, i.e., $\mathcal{D}_{x}$ and $\partial_{i, x}$. The fundamental solution of this first-order elliptic operator is 
\begin{equation*}
	E(x) = \dfrac{1}{\sigma_{n}}\dfrac{\overline{x}}{\arrowvert x\arrowvert^{n + 1}},
\end{equation*} 
where $\sigma_{n}$ is the area of the unit sphere  in $\mathbb{R}^{n + 1}$. Let $G$ be a domain in $\mathbb{R}^{n + 1}$ and $u \in C^{1}(G)$ we say that $u$ is a left-monogenic function (respectively right-monogenic) if $\mathcal{D} u = 0$ (respectively $u\mathcal{D} = 0$) in $G$.

\subsection{Hausdorff and Minkowski Dimensions}\label{SubSecFracDim}
In order to deal with domains with fractal boundaries, we should refresh some basic notions about fractal dimensions. The books \cite{Fal, Mandelbrot, Matt} are recommended as references on this topic. We shall present the notions of Minkowski and Hausdorff dimensions, which are essential tools in this theory. Here we restrict ourselves to the definition of the upper Minkowski dimension. 

\begin{definition}\label{DefMinkDim}$\textbf{(Upper Minkowski dimension)}$
	Let $\textbf{E}$ be a non-empty bounded subset of $\mathbb{R}^{n+1}$ and let $N_{\delta}(\textbf{E})$ be the smallest number of sets of diameter at most $\delta$, covering $\textbf{E}$. The upper Minkowski dimension of $\textbf{E}$ is defined as 
	\begin{equation*}
		\overline{\dim}_{M}\textbf{E} := \limsup_{\delta \rightarrow 0}\dfrac{\log N_{\delta}(\textbf{E}) }{-\log \delta}. 
	\end{equation*}
\end{definition}
The upper Minkowski dimension can also easily be seen to be determined with cubes in a grid, see \cite{Fal}. Suppose $\mathcal{M}_{0}$ denotes a grid covering $\mathbb{R}^{n+1}$ consisting of $(n+1)$-dimensional cubes with sides of length one and vertices with integer coordinates. The grid $\mathcal{M}_{k}$ is obtained from $\mathcal{M}_{0}$ by dividing each of the cubes in $\mathcal{M}_{0}$ into $2^{(n + 1)k}$ different cubes with side lengths $2^{-k}$. Denote by $N_{k}(\textbf{E})$ the number of cubes of the grid $\mathcal{M}_{k}$ which intersect $\textbf{E}$. Then
\begin{equation}\label{EqDefMinkDim}
	\overline{\dim}_{M}\textbf{E} = \limsup_{k \rightarrow \infty}\dfrac{\log N_{k}(\textbf{E}) }{k\log (2)}. 
\end{equation}

Now the concept of the Hausdorff dimension will be introduced. To do that, we need some previous definitions.
\begin{definition}\label{DefHausDim}
	Let $\textbf{E}$ be an arbitrary non-empty set in $\mathbb{R}^{n+1}$. For any $\delta > 0$ and $s \geq 0$, $\mathcal{H}_{\delta}^{s}(\textbf{E})$  is defined as, 
	\begin{equation*}
		\mathcal{H}_{\delta}^{s}(\textbf{E}) := \inf\{\sum_{i = 1}^{\infty}\arrowvert U_{i}\arrowvert^{s}: \{U_{i}\} \ is \ a  \ \delta-covering \ of \ \textbf{E}\},
	\end{equation*} 
	where the infimum is taken over all countable $\delta$-coverings ${U_{i}}$ of $\textbf{E}$ with open or closed balls. We write
	\begin{equation*}
		\mathcal{H}^{s}(\textbf{E}) := \lim _{\delta\rightarrow 0}\mathcal{H}_{\delta}^{s}(\textbf{E}).
	\end{equation*} 
	The Hausdorff dimension of \textbf{E} is defined as
	\begin{equation*}
		\dim _{H}\textbf{E} := \inf\{s \geq 0: \mathcal{H}^{s}(\textbf{E}) = 0\} = \sup\{s \geq 0: \mathcal{H}^{s}(\textbf{E}) = \infty\}.
	\end{equation*}   
\end{definition}
When $s = n+1$ there is a relation between the $(n+1)$-dimensional Lebesgue and Hausdorff measure as we can see in the next theorem. See for instance \cite[pp 28]{Fal}.
\begin{theorem}\label{HausLebeRelation}
	If $\textbf{E}\subset \mathbb{R}^{n+1}$ is a Borel set, then	
	\begin{equation*} 
		\mathcal{H}^{n+1}(\textbf{E}) = \dfrac{1}{\rho_{n+1}} \mathcal{L}^{n+1}(\textbf{E}),
	\end{equation*}
	where $\rho_{n+1}$ is the volume of a $(n+1)$-dimensional ball of diameter one.
\end{theorem} 
In \cite[pp 77]{Matt} is given the next theorem relating the Hausdorff and Minkowski dimensions.
\begin{theorem}
	For the bounded set $\textbf{E} \subset \mathbb{R}^{n+1}$ with topological dimension $n$, we have
	\begin{equation*}
		n \leq \dim_{H}\textbf{E} \leq \overline{\dim}_{M}\textbf{E} \leq n + 1.
	\end{equation*}	
\end{theorem}

\begin{definition}\label{DefFractMan}
	If an arbitrary set $\textbf{E} \subset \mathbb{R}^{n+1}$ with topological dimension $n$ has $\dim_{H}\textbf{E} > n$,  then $\textbf{E}$ is called a fractal set in the sense of Mandelbrot.
\end{definition}
From Definitions \ref{DefHausDim} and \ref{DefFractMan}, we know that a fractal set in the sense of Mandelbrot $\textbf{E}$ satisfies that $\mathcal{H}^{n}(\textbf{E}) = \infty$. Besides, we should note that a bounded set $\textbf{E}$ with $\dim_{H}\textbf{E} = n$ can have $\mathcal{H}^{n}(\textbf{E}) = \infty$. However, classical methods cannot be applied to this kind of set. The ideas developed in this paper are intended to deal with these sets and fractals from Definition \ref{DefFractMan}.

\subsection{Teodorescu Transform and Whitney Extension}
In starting this subsection, we take up some basic properties of the Teodorescu Transform, which will play an essential role in the method developed below, see \cite{GHS08} for more details.

\begin{definition}
	Let $G\subset \mathbb{R}^{n+1}$ be a domain and let $u \in C^{1}(\overline{G})$, the operator defined by $T_{G}$
	\begin{eqnarray*}
		(T_{G}u)(x) = -\int\limits_{G}E(y-x)u(y)dV(y), & x \in \mathbb{R}^{n + 1}, 
	\end{eqnarray*}
	where $dV(y)$ is the volume element, is called the Teodorescu transform. 
\end{definition}  
The next theorem gives us sufficient conditions for the H\"older continuity of the Teodorescu Transform.

\begin{theorem}\label{PropTeo1}
	For $p > n + 1$ and $G$ a domain in $\mathbb{R}^{n+1}$, let $u \in L^{p}(G)$ then\\
	(i) The integral $(T_{G}u)(x)$ exists in the entire $\mathbb{R}^{n + 1}$ and tends to zero for $\arrowvert x\arrowvert \rightarrow \infty$. Besides, $T_{G}u$ is a monogenic function in $\mathbb{R}^{n + 1}\setminus \overline{G}$. Additionally, for bounded domain $G$, we get
	\begin{equation*}
		\|T_{G}u \|_{p} \leq C_{1}(G, p, n)\|u\|_{p}.
	\end{equation*}
	(ii) For $x, y \in \mathbb{R}^{n + 1}$, and $x \neq y$, we have the inequality 
	\begin{equation*}
		\arrowvert(T_{G}u)(x) - (T_{G}u)(y)\arrowvert \leq C_{2}(G, p, n)\|u\|_{p}\arrowvert x - y\arrowvert^{\frac{p - n - 1}{p}}.
	\end{equation*}
\end{theorem}
The following theorem provides conditions for the derivability of the operator $T_{G}u$ over the domain $G$.
\begin{theorem}\label{PropTeo2}
	Let $G$ be a domain and let $u$ be a continuously differentiable function in $G$. Then $T_{G}u$ is also a differentiable function for every $x \in G$ with 
	\begin{equation*}
		\partial_{i}(T_{G}u)(x) = -\int\limits_{G}\partial_{i, x}[E(y - x)]u(y)dV(y) + \overline{e_{i}}\dfrac{u(x)}{n + 1}.
	\end{equation*}
	Particularly, we have the identity 
	\begin{eqnarray*}
		\mathcal{D}(T_{G}u)(x) = u(x),\ \ x \in G.
	\end{eqnarray*}
\end{theorem}

In \cite{ABT2007} can be found the Whitney extension theorem for Clifford valued functions. It is based on the result \cite[pp 174]{St} for real-valued functions, which was stated originally by H. Whitney. This result has enormous importance in this research.

\begin{theorem}\label{ExWhitney}$\textbf{(Whitney extension theorem)}$
	Let $\textbf{E} \subset \mathbb{R}^{n+1}$ be a compact set and let $u \in C^{0, \nu}(\textbf{E})$, with $0 < \nu \leq 1$. Then there exists a function $\widetilde{u} \in C^{0, \nu}(\mathbb{R}^{n+1})$, named Whitney extension operator of $u$, that satisfies\\
	(i) $\widetilde{u}\arrowvert_{E} = u$,\\
	(ii) $\widetilde{u} \in C^{\infty}(\mathbb{R}^{n+1}\setminus \textbf{E})$,\\
	(iii) $\arrowvert \mathcal{D}\widetilde{u}(x) \arrowvert \leq C\dist(x, \textbf{E})^{\nu - 1}$ for $x \in \mathbb{R}^{n+1}\setminus \textbf{E}$.
\end{theorem}

The following theorem is a corollary of a more general result called the Dolzhenko theorem. For the proof, we refer the reader to \cite{APB2007}.
\begin{theorem}\label{CorrDolz}
	Let $G$ be a domain in $\mathbb{R}^{n+1}$ and $\textbf{E} \subset G$ be a compact set. Let be $\mathcal{H}^{n + \mu}(\textbf{E}) = 0$ where $0 < \mu \leq 1$. If $u \in C^{0, \mu}(G)$, and it is monogenic in $G\setminus\textbf{E}$, then $u$ is also monogenic in $G$.   
\end{theorem}

\section{Marcinkiewicz Exponent}\label{Marcinkiewicz}
From now on, let $\mathcal{S}$ be a topologically compact surface, which is the boundary of a Jordan domain in $\mathbb{R}^{n +1}$ that divides it into two domains, the bounded component $G^{+}$ and the unbounded component $G^{-}$ respectively. Let $D \subset \mathbb{R}^{n + 1}$ be a bounded set, which does not intersect the surface $\mathcal{S}$, fractal in general. We define the integral
\begin{equation*}
	I_{p}(D) = \int\limits_{D}\dfrac{dV(x)}{\dist^{p}(x, \mathcal{S})}.
\end{equation*}
When $p = 0$, this integral is the volume of $D$. However, when $p$ is large enough, the integral could diverge.\\
We define the domain $G^{*} := G^{-}\bigcap\{x: \arrowvert x\arrowvert < r\}$, where $r$ is selected such that $\mathcal{S}$ is wholly contained inside the ball of radius $r$. The inner and outer Marcinkiewicz exponent are defined as follows.
\begin{definition}\label{DefMarcExp}
	Let $\mathcal{S}$ be a topologically compact surface which is the boundary of a Jordan domain in $\mathbb{R}^{n +1}$. We define the inner and outer Marcinkiewicz exponent of $\mathcal{S}$, respectively, as
	\begin{equation*}
		\begin{array}{cc}
			\mathfrak{m}^{+}(\mathcal{S}) = \sup \{p: I_{p}(G^{+}) < \infty\},  & \mathfrak{m}^{-}(\mathcal{S}) = \sup\{p: I_{p}(G^{*}) < \infty\}, 
		\end{array} 
	\end{equation*}
	and the (absolute) Marcinkiewicz exponent of $\mathcal{S}$ as,
	\begin{equation*}
		\mathfrak{m}(\mathcal{S}) = \max\{\mathfrak{m}^{+}(\mathcal{S}), \mathfrak{m}^{-}(\mathcal{S})\}.
	\end{equation*} 
\end{definition}

The following lemma plays a significant role in proving the relationship between the Minkowski dimension and the Marcinkiewicz exponent. Here we shall use the Whitney extension decomposition, see \cite{St}.               
\begin{lemma}\label{LemWDescomp}
	Let $\mathcal{S}$ be a topologically compact surface which is the boundary of a Jordan domain in $\mathbb{R}^{n +1}$. Let $w_{k}$ be the number of cubes with edges equal to $2^{-k}$ in the Whitney extension decomposition of $\mathbb{R}^{n +1}\setminus\mathcal{S}$, then
	\begin{equation*}
		w_{k} \leq C2^{kd}
	\end{equation*}
	for each $k \geq m_{d}$ for $m_{d}$ large enough, where $d \in (\overline{\dim}_{M}(\mathcal{S}), n + 1]$, and $C$ is a constant that only depends on $n + 1$.
\end{lemma}
\begin{proof}
	Denote by $w_{k}$ the number of cubes in the grid $\mathcal{M}_{k}$, appearing in the Whitney extension decomposition $\mathcal{F}$ (see \cite{St}). We need to remember that
	\begin{equation}\label{WhitExDes}
		\mathcal{F} = \bigcup_{k}\{Q \in \mathcal{M}_{k}: \,\, Q\cap \Omega_{k} \neq 0, \,\, Q \,\,\,\, is\,\,\,\, maximal\},
	\end{equation}
	where $\Omega_{k}$ is defined as follows
	\begin{equation*}
		\Omega_{k} = \{x: C_{0}2^{-k} \leq \dist(x, \mathcal{S}) \leq C_{0}2^{-k + 1} \},
	\end{equation*}
	here $C_{0} = 2\sqrt{n + 1}$.
	Let $m_{k}(\Omega_{k})$ be the number of cubes of the grid $\mathcal{M}_{k}$, which intersect $\Omega_{k}$, and thus $w_{k} \leq m_{k}(\Omega_{k})$. Suppose that $x \in \Omega_{k}$, then we can find a point $x' \in \mathcal{S}$ separated from $x$ by a distance not greater than $C_{0}2^{-k+1}$. \\
	If $Q$ is a cube of $\mathcal{M}_{k}$, containing $x$, and $Q'$ is a cube of the same grid containing $x'$. Then $Q$ intersects a sphere of radius $C_{0}2^{-k + 1}$ with the center in $Q'$.
	Cubes of $\mathcal{M}_{k}$ intersecting with such spheres lie inside the cube $\widetilde{Q}'$ with edges large enough. Indeed, if we take the cube $Q'' = (1 + C_{0}2^{-k + 1})[Q' - y'] + y'$, where $y'$ is the center of $Q'$, we obtain a cube $Q''$ thicker than $Q'$ by $C_{0}2^{-k+1}$. Hence $Q''$ contains all the spheres centered in a point $x'$ in $Q'$ and the radius equal to $4(\sqrt{n + 1})2^{-k}$. \\
	Let us notice that when $x'$ is in the boundary of $Q'$, the ball with the center in $x'$, and the radius equal to $C_{0}2^{-k+1}$, touches the boundary of $Q''$. Therefore, we need to make $\widetilde{Q}'$ a bit thicker than $Q''$ in order to get all the balls completely contained in $\widetilde{Q}'$. See Figure 1.\\ %\ref{FigLem}.\\
	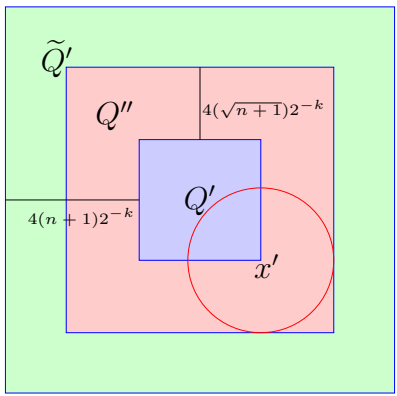
\begin{figure}[h]\label{FigLem}
		\begin{center}
			\begin{tikzpicture}[scale=0.8,draw=blue,fill=blue!20,very thin]	 
				\filldraw[fill=green!20] (3.2,3.2)--(-3.2,3.2)--(-3.2,-3.2)--(3.2,-3.2)--cycle;
				\filldraw[fill=red!20] (2.2,2.2)--(-2.2,2.2)--(-2.2,-2.2)--(2.2,-2.2)--cycle;
				\filldraw (1,1)--(-1,1)--(-1,-1)--(1,-1)--cycle;
				\draw [red] (1,-1) circle (1.2);
				\node at (0,0) {$Q'$};
				\node at (1.1,-1.1) {$x'$};
				\node at (-1.4,1.4) {$Q''$};
				\node at (-2.35,2.35) {$\widetilde{Q}'$};
				\draw [black] (-1,0)--(-3.2,0);
				\draw [black] (0,1)--(0,2.2);
				\node at (1.05,1.5) {\tiny{$4(\sqrt{n + 1})2^{-k}$}};
				\node at (-1.95,-0.3) {\tiny{$4(n + 1)2^{-k}$}};
			\end{tikzpicture}
		\end{center}
		\caption{$\textbf{Two dimensional representation of the cubes $Q'$, $Q''$, and $\widetilde{Q}'$}$.}
	\end{figure}
	It is convenient to have a value in the form $m2^{-k}$, where $m$ is an integer number, in order to get only complete cubes inside $\widetilde{Q}'$. We can choose $\widetilde{Q}' = (1 + 4(n+1)2^{-k})[Q' - y'] + y'$. Now, let us compute the length of the edges of $\widetilde{Q}'$. It is the side of $Q'$ plus twice $4(n+1)2^{-k}$ because it is expanded in both directions, i.e. $[8(n+1) + 1]2^{-k}$.\\
	Consequently, cubes of $\mathcal{M}_{k}$ intersecting with a sphere of the radius equal to $C_{0}2^{-k + 1}$ with the center in $Q'$ stay inside the cube $\widetilde{Q}'$ with edges of length $[8(n+1) + 1]2^{-k}$, and the center coincides with the center of $Q'$. This fact means that for every cube of $\mathcal{M}_{k}$ intersecting $\mathcal{S}$, there exists at most $[8(n+1) + 1]^{n+1}$ cubes of $\mathcal{M}_{k}$ intersecting $\Omega_{k}$. Then
	\begin{equation*}
		w_{k} \leq m_{k}(\Omega_{k}) \leq [8(n+1) + 1]^{n+1}m_{k}(\mathcal{S}).
	\end{equation*}
	From (2) we get that there exists a $N_{0}$ such that for all $k > N_{0}$ we have
	\begin{equation*}
		2^{dk} > m_{k}(\mathcal{S}),
	\end{equation*}
	where $d \in (\overline{\dim}_{M}(\mathcal{S}), n + 1]$ is fixed.
	Consequently,
	\begin{equation*}
		w_{k} \leq m_{k}(\Omega_{k})\leq [8(n+1) + 1]^{n+1}m_{k}(\mathcal{S}) < C2^{dk},
	\end{equation*}
	where $C = [8(n+1) + 1]^{n+1}$.
\end{proof}
In \cite[Lemma 1]{2DKats16}, it is shown, using other tools,  a more general result which particularly implies the next theorem when we restrict ourselves to Lebesgue measure over $\mathbb{R}^{n + 1}$. Here it is shown in a direct way using Lemma \ref{LemWDescomp}.
\begin{theorem}\label{TheoIneqRm}
	Let $\mathcal{S}$ be a topologically compact surface which is the boundary of a Jordan domain in $\mathbb{R}^{n +1}$, then $\mathfrak{m}(\mathcal{S}) \geq n + 1 - \overline{\dim}_{M}(\mathcal{S})$.
\end{theorem}
\begin{proof}
	Let us consider the Whitney extension decomposition (\ref{WhitExDes}) again. We know from \cite{St} that the cubes Q satisfy the inequality 
	\begin{equation}\label{InqCubWhitDes}
		\diam(Q) \leq \dist(Q, \mathcal{S}) \leq 4\diam(Q).
	\end{equation}
	These cubes have edges with lengths equal to $2^{-k}$ where $k \in \mathbb{Z}$ in general. For a fixed cube $Q$ with edge $2^{-k}$ in this decomposition, we infer, from (\ref{InqCubWhitDes}) and since $\diam(Q) = \sqrt{n + 1}\cdotp 2^{-k}$,  that
	\begin{equation*}
		\dfrac{1}{\dist^{p}(x,\mathcal{S})} \leq \dfrac{1}{[\diam(Q)]^{p}} < \dfrac{1}{2^{-kp}}.
	\end{equation*}
	Hence
	\begin{equation*}
		\int_{Q}\dfrac{dV}{\dist^{p}(x,\mathcal{S})} < 2^{k[p - (n + 1)]}.
	\end{equation*}
	We define the values $w^{'}_{k}$ as follows
	
	\begin{equation*}
		w^{'}_{k} := \left\lbrace \begin{array}{cccccccc}
			w_{k}, & if \,\,\, \exists Q_{k} \in \mathcal{M}_{k} & such & that & Q_{k}\cap G^{+}\neq\emptyset, \\
			0, & another \,\,\, case,
		\end{array} \right. 
	\end{equation*}
	where $w_{k}$ is the number of cubes with edges equal to $2^{-k}$ in the Whitney extension decomposition.\\
	Then we have
	\begin{equation*}
		\int\limits_{G^{+}}\dfrac{dV}{\dist^{p}(x,\mathcal{S})} \leq \sum\limits_{Q \in \mathcal{F},\,\,\, Q\cap G^{+}\neq \emptyset }\int\limits_{Q}\dfrac{dV}{\dist^{p}(x,\mathcal{S})} \leq \sum\limits_{k =-\infty}^{\infty}w^{'}_{k}\int\limits_{Q}\dfrac{dV}{\mathrm{\dist}^{p}(x,\mathcal{S})}. 
	\end{equation*}
	However, there is only a finite amount of cubes with edges of length $2^{-k}$  such that $k \leq 0$. Indeed, if $k \leq 0$, then $2^{-k} \geq 1$, and if there are infinitely many cubes with an edge more than or equal to 1, then  the $(n+1)$-dimensional Lebesgue measure of $G^{+}$ would be infinite. In contradiction with the fact that $G^{+}$ is a bounded set in $\mathbb{R}^{n+1}$.\\
	Therefore,
	\begin{equation*}
		\int\limits_{G^{+}}\dfrac{dV}{\dist^{p}(x,\mathcal{S})} \leq C + \sum\limits_{k =1}^{\infty}w^{'}_{k}\int\limits_{Q}\dfrac{dV}{\mathrm{\dist}^{p}(x,\mathcal{S})} < C + \sum\limits_{k = 1}^{\infty}w_{k}2^{k[p - (n +1)]}. 
	\end{equation*} 
	Let $d \in (\overline{\dim}_{M}(\mathcal{S}), n + 1]$, and then from Lemma \ref{LemWDescomp}, we have that 
	\begin{equation*}
		w_{k} \leq B2^{kd},
	\end{equation*}
	for all $k \geq m_{d}$, with $m_{d}$ large enough, and the constant $B$ only depends on $n + 1$. Hence we have
	\begin{equation*}
		\sum\limits_{k = m_{d}}^{\infty}w_{k}2^{k[p - (n +1)]} \leq B\sum\limits_{k = m_{d}}^{\infty}2^{k[p - (n +1) + d]}.
	\end{equation*} 
	Therefore, if the series on the right hand converges, the series on the left side converges. That occurs when is fulfilled the condition
	\begin{equation*}
		p  < (n +1) - d < (n +1) - \overline{\dim}_{M}(\mathcal{S}).
	\end{equation*}
	Consequently,
	\begin{equation*}
		(n +1) - \overline{\dim}_{M}(\mathcal{S}) \leq \mathfrak{m}^{+}(\mathcal{S}).
	\end{equation*}
	An analogous analysis can be done with $G^{*}$ and $\mathfrak{m}^{-}(\mathcal{S})$. 
\end{proof}

\section{A Class of Surfaces in Three Dimensions}\label{Example}
In this section, we construct a class of surfaces in $\mathbb{R}^{3}$. For every possible value of the Minkowski dimension in (2, 3), it is shown that there is a non-numerable class of surfaces with that dimension and such that inequality in Theorem \ref{TheoIneqRm} is strict.
\begin{theorem}\label{TheoInEx}
	Let $\alpha \geq 1$ and $\beta \geq 2$. For each value $d \in (2, 3)$, there exists a non-numerable class of topologically compact surfaces $\mathcal{S}_{\alpha, \beta}$, which are the boundary of a Jordan domain in $\mathbb{R}^{3}$ such that $d = \overline{\dim}_{M}(\mathcal{S}_{\alpha, \beta})$ and $\mathfrak{m}(\mathcal{S}_{\alpha, \beta}) > 3 - d$ for suitable values of $\alpha$ and $\beta$.
\end{theorem}
This construction is similar in spirit to a two-dimensional curve developed in \cite{1DKats16}. That idea on the complex plane goes back at least as far as \cite{BKats83}. The construction follows the simple idea of adding infinitely many three-dimensional rectangles with suitable dimensions to a three-dimensional cube. This begins with a cube $Q = [0, 1]\times[-1, 0]\times[-1, 0]$. Let us fix $\alpha \geq 1$ and $\beta \geq 2$. First, we look at the segment $[0,1]$ in the $x_{1}$ axis. We divide it into infinitely many segments of the form $[2^{-n},2^{-n+1}]$ for each $n \in \mathbb{N}$. Then, for each $n \in \mathbb{N}$, we divide the segments $[2^{-n},2^{-n+1}]$  into $2^{[n\beta]}$ equally spaced segments where $[n\beta]$ is the integer part of $n\beta$. We denote by $x_{nj}$, where $j = 1, 2, ..., 2^{[n\beta]}$, the points determined at the right side of these segments. See Figure 2. \\
%~\ref{lsa}.\\

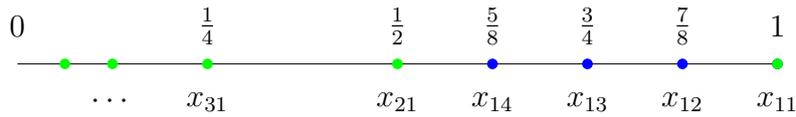
\begin{figure}[h] \label{lsa}
	\usetikzlibrary {calc}
	\begin{center}
		\begin{tikzpicture}[scale=10,fill=blue!20,very thin]		
			\draw (0,0)--(1,0);
			\foreach \j in {1,0.875,0.75,0.625}
			\fill [blue] ($\j*(1,0)$) circle (0.2pt);
			\foreach \i in {1,0.5,0.25,0.125,0.0625}
			\fill [green] ($\i*(1,0)$) circle (0.2pt);
			\node at (1,-0.05) {$x_{11}$};
			\node at (0.5,-0.05) {$x_{21}$};
			\node at (0.25,-0.05) {$x_{31}$};	
			\node at (0.125,-0.05) {$\cdots$};	
			\node at (0.875,-0.05) {$x_{12}$};	
			\node at (0.75,-0.05) {$x_{13}$};
			\node at (0.625,-0.05) {$x_{14}$};	
			\node at (0,0.05) {$0$};
			\node at (1,0.05) {$1$};	
			\node at (0.5,0.05) {$\frac{1}{2}$};	
			\node at (0.25,0.05) {$\frac{1}{4}$};	
			\node at (0.75,0.05) {$\frac{3}{4}$};
			\node at (0.625,0.05) {$\frac{5}{8}$};	
			\node at (0.875,0.05) {$\frac{7}{8}$};	
		\end{tikzpicture}
	\end{center}
	\caption{$\textbf{Distribution of some $x_{nj}$ in the $x_{1}$ axes for $\beta = 2.1$.}$}
\end{figure}
\begin{figure}[h]
	\centering
	\includegraphics[width=0.85\linewidth, height=0.3\textheight]{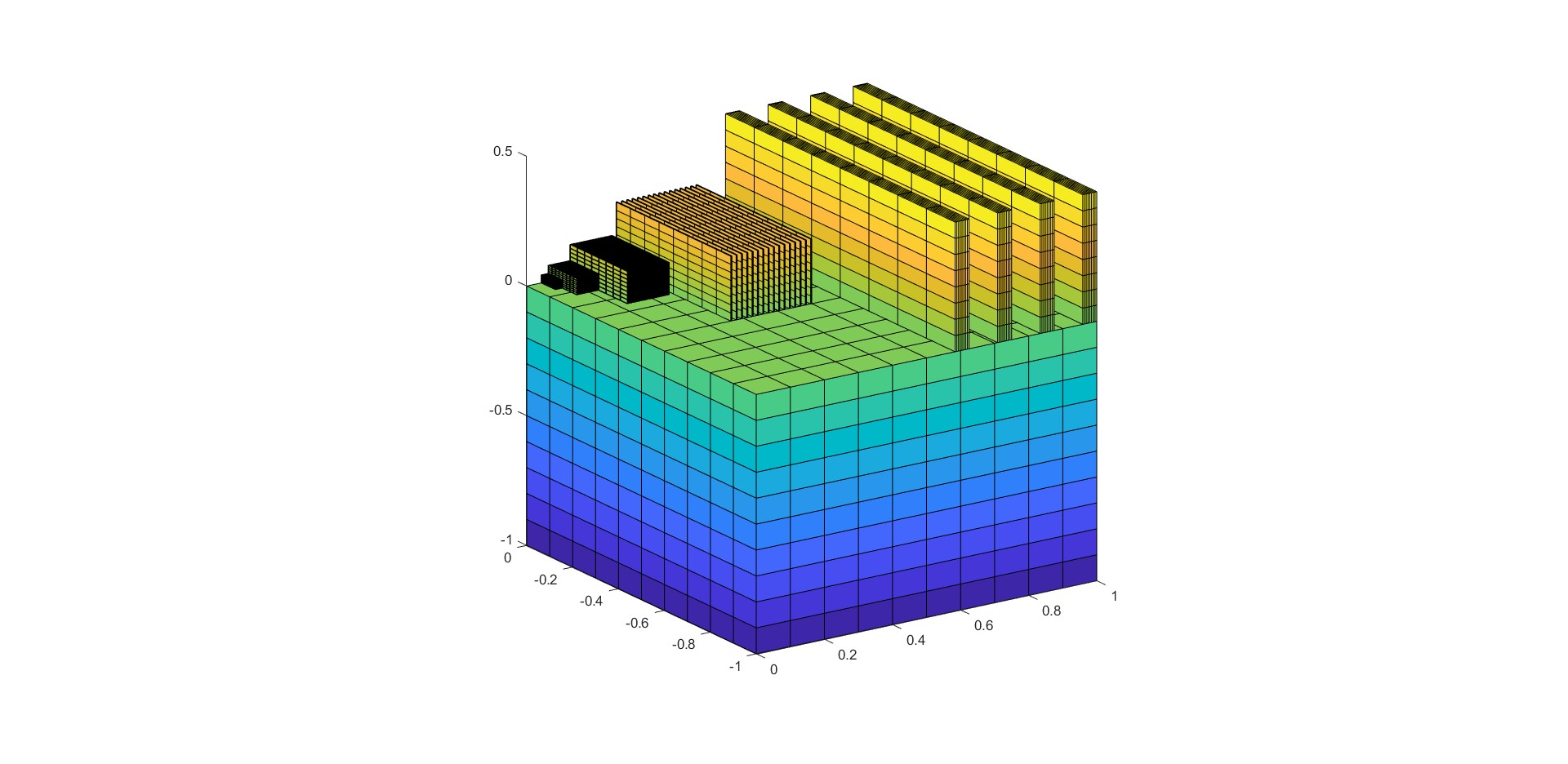}
	\caption[Figure 1. Surface]{The surface $\mathcal{S}_{\alpha, \beta}$ for  $\beta = 2.1$, $\alpha = 1.3$}
	\label{fig:superficie}
\end{figure}

Let $a_{n}$ be the distance between $x_{nj}$ and $x_{n(j+1)}$, i.e. $a_{n} = 2^{-n - [n\beta]}$ and $C_{n} =\frac{1}{2}a_{n}^{\alpha}$. Then let $R_{nj}$ be the following three-dimensional rectangles:
\begin{equation*}
	R_{nj} =  [x_{nj} - C_{n}, x_{nj}]\times[ -2^{-n + 1}, 0]\times[0, 2^{-n}],
\end{equation*}
We define the set 
\begin{equation*}
	T_{\alpha, \beta} = Q\bigcup(\bigcup_{n = 1}^{\infty}\bigcup_{j = 1}^{2^{[n\beta]}}R_{nj}).	
\end{equation*}
We take the surface $\mathcal{S}_{\alpha, \beta} = \partial T_{\alpha, \beta}$. See Figure \ref{fig:superficie}, which was generated using MATLAB, as an illustration. We should note here that the parameter $\beta$ only affects the number of rectangles $R_{nj}$ for each $n \in \mathbb{N}$, while $\alpha$ only affects the width of the rectangles $R_{nj}$.   

\subsection{Minkowski Dimension of the Surfaces $\mathcal{S}_{\alpha, \beta}$}
\begin{figure}[h]
	\centering
	\includegraphics[width=0.85\linewidth, height=0.3\textheight]{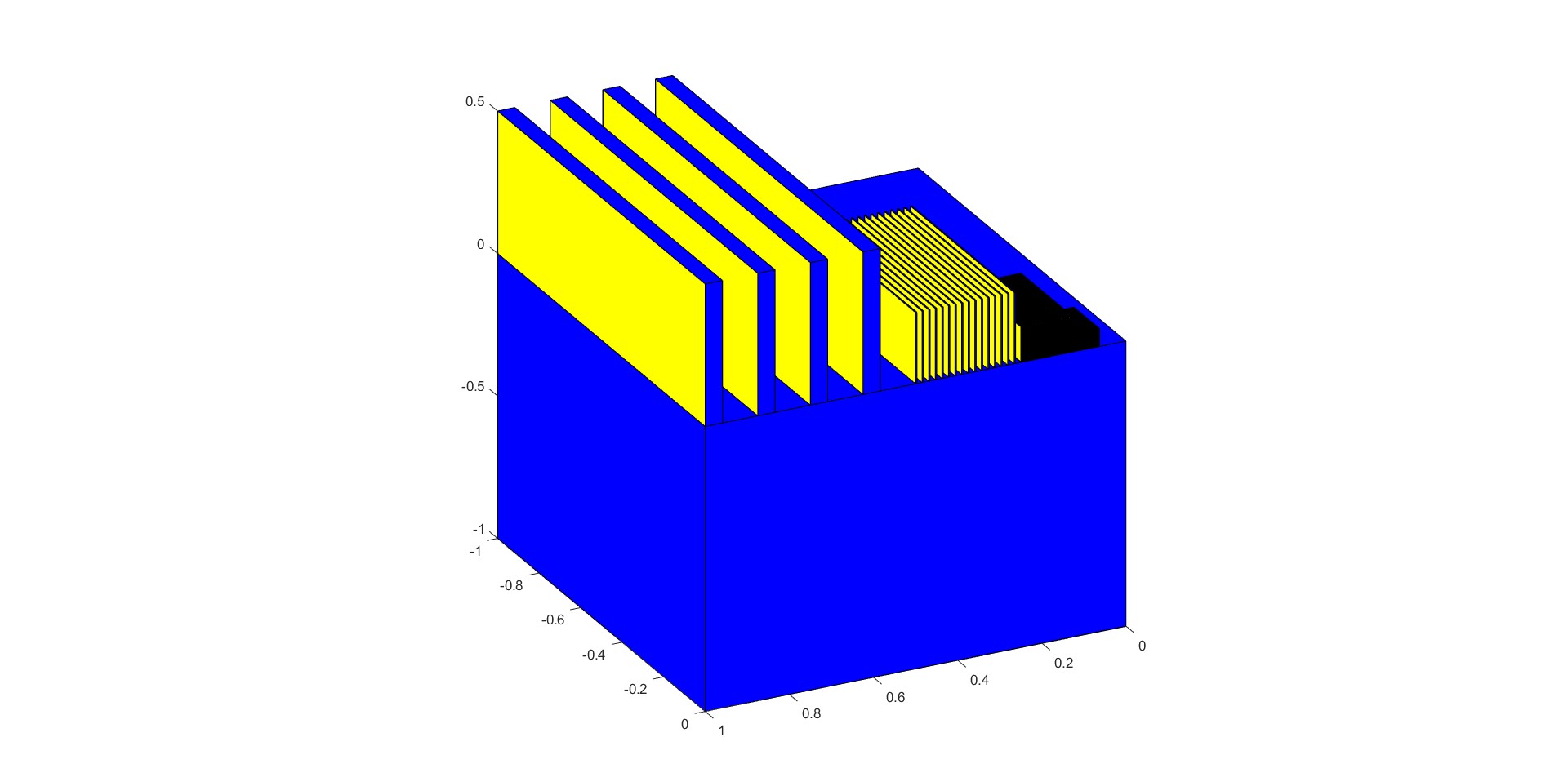}
	\caption[Figure 1. Surface]{Set $A_{\beta}$ into the surface $\mathcal{S}_{\alpha,\beta}$ for $\beta = 2.1$, $\alpha =1.3$}
	\label{fig:Set A into S}
\end{figure}

Now, let us compute the Minkowski dimension of $\mathcal{S}_{\alpha, \beta}$. In order to do that, we shall use the grid $\mathcal{M}_{k}$ defined in \ref{SubSecFracDim}. Many straightforward steps are omitted in order to reduce the exposition. We need to find a lower and an upper bound such that they are equal. To calculate the lower bound, we shall construct a set $A_{\beta}$ such that $A_{\beta} \subset \mathcal{S}_{\alpha, \beta}$ and therefore $ \overline{\dim}_{M}(A_{\beta}) \leq \overline{\dim}_{M}(\mathcal{S}_{\alpha, \beta})$.\\
Let $P_{nj}$ be the two-dimensional rectangles defined as:
\begin{equation*}
	P_{nj} = \{x_{nj}\}\times[-2^{-n+1}, 0]\times[0, 2^{-n}],
\end{equation*}

and the set $A_{\beta}$ is defined as the union 
\begin{equation*}
	A_{\beta} = \bigcup_{n = 1}^{\infty}\bigcup_{j = 1}^{2^{[n\beta]}}P_{nj},
\end{equation*}
see Figure \ref{fig:Set A into S}. This was created with the software MATLAB.  

From construction, we know that $A_{\beta} \subset \mathcal{S}_{\alpha, \beta}$. We are going to find a lower bound for $\overline{\dim}_{M}(A_{\beta})$. To do that, let us focus on the distance between $P_{nj}$ and $P_{nj+1}$. It is equal to $a_{n} = 2^{-n -[n\beta]}$. If $k > n$, and $a_{n} > 2^{-k}$, a cube cannot touch two of these rectangles. The quantity of cubes in $\mathcal{M}_{k}$ that cover a single two-dimensional rectangle $P_{nj}$ is $2(\frac{2^{-n}}{2^{-k}})^{2}$ because these two-dimensional rectangles have lengths of $2^{-n +1}$ and widths of $2^{-n}$.\\
There are $2^{[n\beta]}$ rectangles $P_{nj}$ for a fixed $n$. Therefore $2^{[n\beta] +1}(\frac{2^{-n}}{2^{-k}})^{2}$ cubes are needed to cover all the $P_{nj}$ for a fixed $n$. Then we have
\begin{equation*}
	N_{k}(A_{\beta}) \geq 2\cdotp\sum\limits_{a_{n}> 2^{-k},\,\,\, k>n}2^{[n\beta]+2k-2n},
\end{equation*}
where $N_{k}(A_{\beta})$ is the minimal number of cubes of the grid $\mathcal{M}_{k}$ which cover $A_{\beta}$.\\
Denote by $B_{k}$ the integer defined by the condition 
\begin{equation}\label{IneBk}
	\dfrac{k}{1+\beta}-1 \leq B_{k} < \dfrac{k}{1+\beta}.
\end{equation} 
It is not difficult to show that the condition $a_{n} > 2^{-k}$ is fulfilled if and only if $n = 1, 2,..., B_{k}$. Next, we get
\begin{equation*}
	\sum\limits_{a_{n}> 2^{-k},\,\,\, k>n}2^{[n\beta]+2k-2n} = 2^{2k}\sum\limits_{n = 1}^{B_{k}}2^{[n\beta]-2n} \geq 2^{2k-1}\sum\limits_{n = 1}^{B_{k}}2^{n(\beta-2)} \geq C2^{\frac{3k\beta}{\beta+1}},
\end{equation*}
where $C$ does not depend on $k$. Therefore
\begin{equation*}
	\overline{\dim}_{M}(S_{\alpha,\beta}) \geq \overline{\dim}_{M}(A_{\beta}) \geq \frac{3\beta}{\beta + 1}.
\end{equation*}
We need to find an accurate upper bound for $\overline{\dim}_{M}(S_{\alpha, \beta})$. In order to do that, we define the sets $\Lambda_{n} := \bigcup\limits_{j = 1}^{2^{[n\beta]}}[\partial R_{nj}\setminus(\partial R_{nj})\arrowvert_{x_{3} = 0}]$ and $\Lambda := \bigcup\limits_{n = 1}^{\infty}\Lambda_{n}$. Defining $\widehat{Q} := \partial Q\setminus[\bigcup\limits_{n = 1}^{\infty}\bigcup\limits_{j = 1}^{2^{[n\beta]}}(\partial R_{nj})\arrowvert_{x_{3} = 0}]$, we can observe that $\mathcal{S}_{\alpha, \beta} = \widehat{Q}\cup\Lambda$. We shall focus first on $\Lambda$.\\
We are going to consider three cases. The first one is when $n \geq k$, the second one will be if $n < k$ and $C_{n} > 2^{-k}$, and the last one is if $n<k$ and $C_{n} \leq 2^{-k}$. From construction, the surfaces $\Lambda_{n}$, with $n > k$, are covered by one cube of the grid $\mathcal{M}_{k}$. While the surface $\Lambda_{k}$ is covered by two of these cubes. As above, if $n < k$ and $C_{n} > 2^{-k}$, then $2^{[n\beta] +2}(\frac{2^{-n}}{2^{-k}})^{2}$ cubes are needed to cover the sides of the $R_{nj}$ parallel to $x_{1} = 0$, in $\Lambda_{n}$.\\
No more than $2(\frac{2^{-n}}{2^{-k}})^{2}$ cubes are needed to cover the two-dimensional rectangles in $\Lambda_{n}$ parallel to the coordinate plane $x_{2} = 0$. In addition, $2(\frac{2^{-n}}{2^{-k}})^{2}$ cubes are enough to cover the two-dimensional rectangles in $\Lambda_{n}$ parallels to $x_{3} = 0$.\\
If $n<k$ and $C_{n} \leq 2^{-k}$, let us analyze two more cases; when $a_{n} - C_{n} \leq 2^{-k}$ and $a_{n} - C_{n} > 2^{-k}$.\\
Following the same idea, we get that if $C_{n} \leq 2^{-k}$, $k > n$, and $a_{n} - C_{n} \leq 2^{-k}$, then $2(\frac{2^{-n}}{2^{-k}})^{3}$ cubes in $\mathcal{M}_{k}$ are enough to cover $\Lambda_{n}$.\\
If $C_{n} \leq 2^{-k}$, $k > n$, and $a_{n} - C_{n} > 2^{-k}$, then no more than
$2^{[n\beta] +2}(\frac{2^{-n}}{2^{-k}})^{2}$ cubes are needed to cover the faces of $R_{nj}$'s parallel to $x_{1} = 0$ in $\Lambda_{n}$. No more than $2(\frac{2^{-n}}{2^{-k}})^{2}$ cubes are needed to cover the two-dimensional rectangles in $\Lambda_{n}$ parallel to the coordinate plane $x_{2} = 0$. Besides,  $(\frac{2^{-n}}{2^{-k}})^{2}$ cubes are enough to cover the two-dimensional rectangles in $\Lambda_{n}$ parallels to $x_{3} = 0$.\\
Finally, with $6(\frac{1}{2^{-k}})^{2}$ cubes of $\mathcal{M}_{k}$, we can cover $\widehat{Q}$. As a consequence, we get
\begin{equation*}
	N_{k}(S_{\alpha,\beta}) \leq 3 + 6\cdotp2^{2k} + 4\sum\limits_{C_{n}> 2^{-k},\,\,\, k>n}2^{[n\beta]+2k-2n} + 4\sum\limits_{C_{n}> 2^{-k},\,\,\, k>n}2^{2k-2n} + 
\end{equation*}

\begin{equation*}
	+ 2\sum\limits_{C_{n}\leq 2^{-k},\,\,\, a_{n} - C_{n} \leq 2^{-k},\,\,\, k>n}2^{3k-3n} + 4\sum\limits_{C_{n}\leq 2^{-k} < a_{n} - C_{n},\,\,\, k>n}2^{[n\beta]+2k-2n} +
\end{equation*}

\begin{equation*}
	+ 4\sum\limits_{C_{n}\leq 2^{-k} < a_{n} - C_{n} \,\,\,	k>n}2^{2k-2n}. 
\end{equation*}

Using the conditions on the sums, it is possible to get estimates greater than those obtained before and then
\begin{equation*}
	N_{k}(S_{\alpha,\beta}) \leq 3 + 6\cdotp2^{2k} + 8\sum\limits_{2^{-k} < a_{n},\,\,\, k>n}2^{[n\beta]+2k-2n} + 8\sum\limits_{2^{-k} < a_{n},\,\,\, k>n}2^{2k-2n} + 
\end{equation*}

\begin{equation*}
	+ 2\sum\limits_{\frac{a_{n}}{2} \leq 2^{-k},\,\,\, k>n}2^{3k-3n}. 
\end{equation*}
Using $B_{k}$ defined in (\ref{IneBk}) for those sums under the conditions $2^{-k} < a_{n}, k>n$; and the integer $H_{k}$ defined
\begin{equation*}
	\dfrac{k-1}{1+\beta}-1 \leq H_{k} < \dfrac{k-1}{1+\beta},
\end{equation*} 
for the sum under the conditions $\frac{a_{n}}{2} \leq 2^{-k}, k>n$, we can obtain through simple estimates the next inequality
\begin{equation*}
	N_{k}(S_{\alpha,\beta}) \leq D(k)2^{\frac{3k\beta}{\beta+1}},
\end{equation*}
where $D(k) = ak + c$; here $a$ and $c$ only depend on $\beta$.
Hence
\begin{equation*}
	\overline{\dim}_{M}(\mathcal{S}_{\alpha, \beta}) \leq \dfrac{3\beta}{\beta+1}.
\end{equation*}
Consequently,
\begin{equation*}
	\overline{\dim}_{M}(\mathcal{S}_{\alpha, \beta}) = \dfrac{3\beta}{\beta+1}.
\end{equation*}

\subsection{Marcinkiewicz Exponent of the Surfaces $\mathcal{S}_{\alpha, \beta}$}
Here we shall compute the Marcinkiewicz exponent. Again many straightforward steps are omitted to shorten the exposition. In order to do that, we divide $T_{\alpha, \beta}$ into regions where we can express the function $\dist(x,\mathcal{S}_{\alpha, \beta})$ in terms of elementary functions. In $Q$, we can draw the planes which bisect the dihedral angle between two adjacent faces of $Q$. All these planes intersect each other at the point $A = (\frac{1}{2}, -\frac{1}{2}, -\frac{1}{2})$. In that way, we divide $Q$ into six different right square pyramids with vertex at $A$.\\ 
We call $G^{+}_{1}$ and $G^{+}_{2}$ to the pyramids with base $\{0\}\times[-1, 0]\times[-1, 0]$ and its parallel face respectively. Similarly, $G^{+}_{3}$ and $G^{+}_{5}$ have bases $ [0, 1]\times\{-1\}\times[-1, 0]$ and $[0, 1]\times[-1, 0]\times\{-1\}$ while the bases of $G^{+}_{4}$ and $G^{+}_{6}$ are its parallel faces, respectively. Finally, let be $G^{+}_{7} = \bigcup\limits_{n = 1}^{\infty}\bigcup\limits_{j = 1}^{2^{[n\beta]}} R_{nj}$. Hence we have $G^{+} =  \bigcup\limits_{i = 1}^{7}G^{+}_{i}$.
Due to the fact that the faces of the right square pyramids bisect the dihedral angles between adjacent faces of the cube $Q$, we get
\begin{equation*}
	\begin{array}{ccc}
		\dist^p(x, \mathcal{S}_{\alpha, \beta})\arrowvert_{G^{+}_{1}} = x_{1}^{p}, & \dist^p(x, \mathcal{S}_{\alpha, \beta})\arrowvert_{G^{+}_{2}} = (1 - x_{1})^{p},\\
		dist^p(x, \mathcal{S}_{\alpha, \beta})\arrowvert_{G^{+}_{4}} = \arrowvert x_{2}\arrowvert^{p}, & \dist^p(x, \mathcal{S}_{\alpha, \beta})\arrowvert_{G^{+}_{3}} = \arrowvert-1 - x_{2}\arrowvert^{p},\\
		\dist^p(x, \mathcal{S}_{\alpha, \beta})\arrowvert_{G^{+}_{6}} \geq \arrowvert x_{3}\arrowvert^{p}, & \dist^p(x, \mathcal{S}_{\alpha, \beta})\arrowvert_{G^{+}_{5}} = \arrowvert-1 - x_{3}\arrowvert^{p}.\\
	\end{array}
\end{equation*}
Since the faces of the pyramids have null volume, we get that
\begin{equation*}
	\int\limits_{G^{+}}\dfrac{dV}{\dist^{p}(x,\mathcal{S}_{\alpha, \beta})} = \sum_{i = 1}^{7} \int\limits_{G^{+}_{i}}\dfrac{dV}{\dist^{p}(x, \mathcal{S}_{\alpha,\beta})}.
\end{equation*}
Furthermore, because $\dist^{p}(x, \mathcal{S}_{\alpha, \beta})$ is a positive function, the integral in the left hand converges if and only if the seven integrals in the sum in the right hand converge.\\
It is possible to show through direct computations that 
\begin{equation*}
	\int\limits_{G^{+}_{1}}\dfrac{dV}{x_{1}^{p}} < \infty,
\end{equation*}  
if and only if $p < 1$. Hence we only need to analyze these values of $p$ in the following integrals. Analogous computations can be done to obtain that the integrals over the regions $G^{+}_{i}$, where $i = 2, ..., 6$ converge when $p < 1$.\\
On the other hand, for the integral over the region $G^{+}_{7}$, we have that
\begin{equation*}
	\int\limits_{G^{+}_{7}}\dfrac{dV}{\dist^{p}(x, \mathcal{S}_{\alpha,\beta})} = \sum_{n = 1}^{\infty}\sum_{j = 1}^{2^{[n\beta]}} \int\limits_{R_{nj}}\dfrac{dV}{\dist^{p}(x, \mathcal{S}_{\alpha,\beta})}.
\end{equation*}
In order to compute the integral over $R_{nj}$, we divide that region in the same way that in the cube $Q$. By drawing the planes that bisect the dihedral angles at the edges, we get regions where the function $\dist^{p}(x, \mathcal{S}_{\alpha,\beta})$ can be represented through elementary functions.\\
After doing the tedious calculations, we can reduce the convergence of the integral over $G^{+}_{7}$ when $p<1$ to the convergence of the series
\begin{equation*}
	\sum_{n = 1}^{\infty}2^{[n\beta] - 2n}\left(\frac{C_{n}}{2}\right)^{1-p},
\end{equation*}
which converges if and only if converges the series
\begin{equation*}
	\sum_{n = 1}^{\infty}2^{n\beta - 2n - (1-p)\alpha(n + n\beta)}.
\end{equation*}
This geometric series converges if and only if the condition  
\begin{equation*}
	p < 1 - \dfrac{\beta - 2}{\alpha(\beta+1)},
\end{equation*}
is fulfilled. Thus, we have that the inner Marcinkiewicz exponent is 
\begin{equation*}
	\mathfrak{m}^{+}(\mathcal{S}_{\alpha, \beta}):= \sup\{p>0: I_{p}(G^{+})<\infty\} = 1 - \dfrac{\beta - 2}{\alpha(\beta+1)}.
\end{equation*}
Also, we obtain that the absolute Marcinkiewicz exponent satisfy
\begin{equation*}
	\mathfrak{m}(\mathcal{S}_{\alpha, \beta}) := \max\{\mathfrak{m}^{+}(\mathcal{S}_{\alpha, \beta}), \mathfrak{m}^{-}(\mathcal{S}_{\alpha, \beta})\} \geq \mathfrak{m}^{+}(\mathcal{S}_{\alpha, \beta}) = 1 - \dfrac{\beta - 2}{\alpha(\beta+1)}.
\end{equation*}

\subsection{Remarks about the Surfaces $\mathcal{S}_{\alpha, \beta}$}
Now we are able to prove Theorem \ref{TheoInEx}.\\
\begin{proof}[Proof of Theorem~{\upshape\ref{TheoInEx}}]
	If $\alpha > 1 $ and $\beta > 2$ then
	\begin{equation*}
		\mathfrak{m}(\mathcal{S}_{\alpha, \beta}) \geq \mathfrak{m}^{+}(\mathcal{S}_{\alpha, \beta}) = 1 - \dfrac{\beta - 2}{\alpha(\beta+1)} > 1 - \dfrac{\beta - 2}{\beta+1} = 3 - \dfrac{3\beta}{\beta + 1} = 3 - \overline{\dim}_{M}(\mathcal{S}_{\alpha, \beta}).
	\end{equation*}
	For each $d \in (2, 3)$, let be $\beta = \frac{d}{3 - d}$ then $\overline{\dim}_{M}(\mathcal{S}_{\alpha, \beta}) = d$ for each $\alpha > 1$, i. e. a non-numerable family. 
\end{proof}
On the other hand, as a trivial conclusion we see that when $\alpha = 1 $ or $\beta = 2$ we have that $\mathfrak{m}^{+}(\mathcal{S}_{\alpha, \beta}) = 3 - \overline{\dim}_{M}(\mathcal{S}_{\alpha, \beta})$.\\
We can also note that when $\beta = 2$ then  $2 \leq \dim_{H}(\mathcal{S}_{\alpha, 2}) \leq \overline{\dim}_{M}(\mathcal{S}_{\alpha, 2}) = 2$ and consequently $\dim_{H}(\mathcal{S}_{\alpha, 2}) = 2$. On the other hand, the 2-Hausdorff measure is $\mathcal{H}^{2}(\mathcal{S}_{\alpha, 2}) = \infty$, because $\mathcal{H}^{2}(\mathcal{S}_{\alpha, 2}) \geq \mathcal{H}^{2}(A_{2})$ and from Theorem \ref{HausLebeRelation} we have that $\mathcal{H}^{2}(A_{2}) = \infty$ . Therefore, $\mathcal{S}_{\alpha, 2}$ is not a fractal in the sense of Mandelbrot. However, classical methods cannot be applied to it, even those developed for non-smooth surfaces.\\
Even though it is impossible to draw a hypersurface like this example in dimensions higher than three, we are able to describe it analytically. Indeed, let $ Q = [0, 1]\times[0, 1]\times[0, 1]\times\cdots\times[-1, 0]$ be a $(n+1)$-dimensional cube. Additionally, let $R_{mj}$ be the $(n+1)$-dimensional rectangles given by
\begin{equation*}
	R_{mj} =  [x_{mj} - C_{m}, x_{mj}]\times[0, 2^{-m}]\times\cdots\times[0, 2^{-m}],
\end{equation*}
a product of $(n+1)$ segments. Then we analogously define 
\begin{equation*}
	T^{n+1}_{\alpha, \beta} = Q\bigcup(\bigcup_{m = 1}^{\infty}\bigcup_{j = 0}^{2^{[m\beta]} - 1}R_{mj}),	
\end{equation*}
where the hypersurface $\mathcal{S}^{n+1}_{\alpha, \beta} = \partial T^{n+1}_{\alpha, \beta}$. We should note that for $n+1 = 3$ this surface is pretty similar to the one in the Figure \ref{fig:superficie}.

\section{Jump Problem in Fractal Domains}\label{JumpProblem} 
Throughout this section, the following temporary notation will be used. Let $\mathcal{S}$ be a topologically compact surface that is the boundary of a Jordan domain $G^{+} \subset \mathbb{R}^{n + 1}$, and let be $G^{-} :=(\mathbb{R}^{n+1}\bigcup\{\infty\})\setminus (G^{+}\bigcup\mathcal{S})$. The Jump Problem in Clifford Analysis is stated as follows: Given a continuous function $f$ defined on $\mathcal{S}$, to find two functions $\Phi^{+}(x)$, monogenic in the domain $G^{+}$, and $\Phi^{-}(x)$ monogenic in the domain $G^{-}$, such that they have continuous limit values in the boundary $\mathcal{S}$ and  there satisfy the relation
\begin{eqnarray}\label{EqJumpProb}
	\Phi^{+}(x) - \Phi^{-}(x) = f(x), & x \in \mathcal{S},
\end{eqnarray}
with $\Phi(\infty) = 0$.
If $\mathcal{S}$ is a fractal surface, it is impossible to use the cliffordian Cauchy type integral to solve the problem (\ref{EqJumpProb}). In the context of Clifford Analysis, we have the following result, which generalizes \cite[Theorem 2]{1DKats16} to higher dimensions.
\begin{theorem}\label{TheoSolvCond}
	Let $\mathcal{S}$ be a topologically compact surface which is the boundary of a Jordan domain in $\mathbb{R}^{n +1}$,  and let $f \in C^{0, \nu}(\mathcal{S})$. If 
	\begin{equation}\label{EqCndMar}
		\nu > 1 - \dfrac{\mathfrak{m}(\mathcal{S})}{n + 1},
	\end{equation}
	then the jump problem (\ref{EqJumpProb}) is solvable.
\end{theorem}
\begin{proof} 
	First, we consider the inner Marcinkiewicz exponent $\mathfrak{m}^{+}(\mathcal{S})$. We look for sufficient conditions such that the Whitney extension $\widetilde{f}$ of $f$ satisfies that $\mathcal{D}\widetilde{f} \in$ L$^{p}(G^{+})$ with $p > n + 1$. 
	Indeed, from Theorem \ref{ExWhitney}, we have
	\begin{equation*}
		\int\limits_{G^{+}}\arrowvert\mathcal{D}\widetilde{f}(x)\arrowvert^{p}dV(x) \leq C\int\limits_{G^{+}}\dfrac{dV(x)}{\dist(x, \mathcal{S})^{p(1 - \nu)}}.
	\end{equation*}
	
	From Definition \ref{DefMarcExp}, we have that the above right-hand integral converges for $p < \frac{\mathfrak{m}^{+}(\mathcal{S})}{1 - \nu}$. Then we need that $n + 1 < \frac{\mathfrak{m}^{+}(\mathcal{S})}{1 - \nu}$, or equivalently 
	\begin{equation*}
		\nu > 1 - \dfrac{\mathfrak{m}^{+}(\mathcal{S})}{n + 1}.
	\end{equation*} 
	
	Note that this is a sufficient condition for $\mathcal{D}\widetilde{f} \in$ L$^{p}(G^{+})$ with $p > n + 1$. Next, let us consider the function  
	\begin{equation}\label{EqSolJmPro}
		\begin{array}{cc}
			\Phi(x) = \widetilde{f}(x)\chi(x) - (T_{G^{+}}\mathcal{D} \widetilde{f})(x), & x \in \mathbb{R}^{n+1},
		\end{array}
	\end{equation}
	where $\chi(x)$ is the characteristic function of $G^{+}$. We shall show that, under condition (\ref{EqCndMar}), function (\ref{EqSolJmPro}) is a solution to the jump problem.\\
	Indeed, we have that
	\begin{equation*}
		\begin{array}{cc}
			\Phi^{-}(x) = - (T_{G^{+}}\mathcal{D} \widetilde{f})(x), & x \in G^{-}.
		\end{array}
	\end{equation*}
	From Theorem \ref{PropTeo1}, we get that $\Phi^{-}(x)$ is a monogenic function over $G^{-}$, vanish at infinity, and also $\Phi^{-}(x) \in C^{0, \alpha}(\overline{G^{-}})$, with $\alpha = \frac{p - n - 1}{p}$. Consequently, $\Phi^{-}$ is a continuous function over $\overline{G^{-}}$.\\ 
	On the other hand,
	\begin{equation*}
		\begin{array}{cc}
			\Phi^{+}(x) = \widetilde{f}(x) - (T_{G^{+}}\mathcal{D} \widetilde{f})(x), & x \in G^{+},
		\end{array}
	\end{equation*}
	from Theorem \ref{PropTeo1} we know that $(T_{G^{+}}\mathcal{D} \widetilde{f})(x) \in C^{0, \alpha}(\overline{G^{+}})$ with $\alpha = \frac{p - n - 1}{p}$. 
	Besides, we know that $\widetilde{f} \in C^{0, \nu}(\mathbb{R}^{n + 1})$ thus $\Phi^{+}(x)$ is a continuous function over $\overline{G^{+}}$.  From Theorem \ref{PropTeo2}, we get $\mathcal{D}\Phi^{+}(x) = 0$ over $G^{+}$. Finally, we can verify directly that the function $\Phi(x)$ satisfies the boundary condition over $\mathcal{S}$.\\  
	On the other hand, for the outer Marcinkiewicz exponent $\mathfrak{m}^{-}(\mathcal{S})$, we suppose that $\mathcal{S}$ is entirely contained inside the ball $K_{1} = \{x: \arrowvert x\arrowvert < r_{1}\}$. Let be $r > r_{1}$, and $K = \{z: \arrowvert x\arrowvert < r\}$. Besides, let $w(x)$ be a function in $C^{\infty}(\mathbb{R}^{n + 1})$ equal to 1 over $K_{1}$, equal to 0 outside of $K$, and  $0 \leq w(x) \leq 1$. Let be $G^{*} = G^{-}\bigcap K$ and $f^{*} = \widetilde{f}w$, we can observe that $\arrowvert\mathcal{D} f^{*}\arrowvert \leq \frac{C}{(\dist(x, \mathcal{S}))^{1 - \nu}}$. In a similar way, we get that under the condition $\nu > 1 - \frac{\mathfrak{m}^{-}(\mathcal{S})}{n + 1}$, the function 
	\begin{equation*}
		\Phi(x) = f^{*}(x)\chi^{*}(x) - (T_{G^{*}}\mathcal{D} f^{*})(x),
	\end{equation*}
	where $\chi^{*}$ is the characteristic function of $G^{*}$, is a solution to the jump problem. 
\end{proof}
From Theorem \ref{TheoIneqRm}, it follows that Theorem \ref{TheoSolvCond} improves the existing conditions for the solvability of the jump problem. Additionally, using Theorem \ref{CorrDolz}, we obtain the following unicity conditions. 
\begin{theorem}\label{TheoUniMar}
	Let $\mathcal{S}$ be a topologically compact surface which is the boundary of a Jordan domain in $\mathbb{R}^{n + 1}$, and let $f \in C^{0, \nu}(\mathcal{S})$, with $\nu > 1 - \dfrac{\mathfrak{m}(\mathcal{S})}{n + 1}$ and 
	\begin{equation}\label{EqUniCndMar}
		\dim_{H}\mathcal{S} - n < \mu < 1 - \dfrac{(n + 1)(1 - \nu)}{\mathfrak{m}(\mathcal{S})}.
	\end{equation}
	Then the solution to the jump problem (\ref{EqJumpProb}) is unique in the classes $C^{0, \mu}(\overline{G^{+}})$ and $C^{0, \mu}(\overline{G^{-}})$.
\end{theorem}
The unicity in Theorem \ref{TheoUniMar} is assumed when there exists a value of $\mu$ such that condition (\ref{EqUniCndMar}) is fulfilled.

\begin{further}
	The results obtained here are also valid in the context of vectorial Clifford analysis. There we need to use the properties of the Teodorescu transform written in the vectorial sense that can be found in \cite{APB2007, GS1997} and analogous reasoning. We should note that in the theorems with this approach, we need to change the dimension of the space from $(n+1)$ to $n$. 
\end{further}

\subsubsection*{Acknowledgments} The author gratefully acknowledges the financial support of the Postgraduate Study Fellowship of the Consejo Nacional de Ciencia y Tecnología (CONACYT) (Grant Number 957110).

\end{document}